\documentclass[12pt]{amsart}
\usepackage{amssymb,amsmath,amsthm,mathtools,amsfonts,times,hyperref,mathrsfs,multirow}
\usepackage{pgfplots}
\pgfplotsset{compat=1.15}
\usepackage{mathrsfs,mathdots}
\usetikzlibrary{arrows}
\usepackage{enumerate}
\usepackage{graphicx}
\usepackage{array}
\usepackage{ytableau}
\usepackage{pstricks,pst-node,graphicx}
\usepackage{tikz,varwidth}
\usepackage{setspace}
\usepackage{float}
\usetikzlibrary{calc}
\usepackage[super]{nth}
\usepackage{tikz}
\usepackage{extarrows}
\usepackage{pgfplots}
\pgfplotsset{compat=1.15}
\usetikzlibrary{arrows}

\newtheorem{theorem}{Theorem}[section]

\newtheorem{cor}[theorem]{Corollary}
\newtheorem{conj}[theorem]{Conjecture}

\theoremstyle{definition}

\theoremstyle{remark}
\newtheorem{remark}{Remark}
\numberwithin{equation}{section}

\allowdisplaybreaks

\newcommand\numberthis{\addtocounter{equation}{1}\tag{\theequation}}
\newcommand\restr[2]{{
  \left.\kern-\nulldelimiterspace 
  #1 
  \littletaller 
  \right|_{#2} 
  }}
\newcommand{\littletaller}{\mathchoice{\vphantom{\big|}}{}{}{}}  

\begin{document}

\title[]
{Extension of Bressoud's generalization of Borwein's conjecture and some exact results}

\author{Alexander Berkovich}
\address{Department of Mathematics, University of Florida, Gainesville
FL 32611, USA}
\email{alexb@ufl.edu}
\author{Aritram Dhar}
\address{Department of Mathematics, University of Florida, Gainesville
FL 32611, USA}
\email{aritramdhar@ufl.edu}

\dedicatory{Dedicated to George E. Andrews and Bruce C. Berndt in celebration of their \nth{85} birthdays}

\date{\today}

\subjclass[2020]{05A15, 05A17, 05A30, 11P81, 11P84}             

\keywords{hook differences, Bressoud's conjecture, positivity-preserving transformations, $q$-series with non-negative coefficients}

\begin{abstract}
In this paper, we conjecture an extension to Bressoud's 1996 generalization of Borwein's famous 1990 conjecture. We then state a few infinite hierarchies of non-negative $q$-series identities which are interesting examples of our proposed conjecture and Bressoud's generalized conjecture. Finally, using certain positivity-preserving transformations for $q$-binomial coefficients, we prove the non-negativity of the infinite families.
\end{abstract}
\maketitle

\section{Introduction}\label{s1}
A \textit{partition} $\pi$ is a non-increasing finite sequence $\pi = (\lambda_1,\lambda_2,\dots)$ of positive integers. The elements $\lambda_i$ that appear in the sequence $\pi$ are called \textit{parts} of $\pi$. The \textit{number of parts} of $\pi$ is denoted by $\#(\pi)$. The sum of all the parts of a partition $\pi$ is called the \textit{size} of this partition and is denoted by $|\pi|$. We say $\pi$ is a partition of $n$ if its size is $n$. The empty sequence $\emptyset$ is considered as the unique partition of zero.\quad\\\par The Young diagram of $\pi$ is a convenient way of representing $\pi$ graphically wherein the parts of $\pi$ are depicted as rows of unit squares which are called cells. Let $\pi$ be a partition whose Young diagram has a node in the $i$-th row and $j$-th column. We call this node the $(i,j)$-\textit{node}. We define the \textit{hook difference} at the $(i,j)$-th node to be the number of nodes in the $i$-th row of $\pi$ minus the number of nodes in the $j$-th column of $\pi$. For example, if $\pi$ is $5+3+1$, then the Young diagram of $\pi$ is
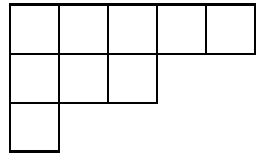
\begin{figure}[H]
\centering
\ydiagram{5,3,1}
\end{figure}
and the hook differences at each node are
\begin{figure}[H]
\centering
\ytableaushort
{{2}{3}{3}{4}{4},{0}{1}{1},{-2}}
* {5,3,1}
\par\quad\\\par\quad\\
\caption{Hook differences of the partition $\pi = (5,3,1)$}
\label{fig1}
\end{figure}
We say that the $(i,j)$-th node lies on diagonal $c$ if $i-j=c$. For example, if $\pi$ is $5+3+1$, then the diagonals are
\vspace{-1cm}
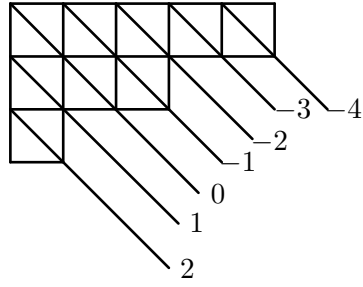
\begin{figure}[H]
\begin{tikzpicture}[line cap=round,line join=round,>=triangle 45,x=0.7cm,y=0.7cm]
\clip(-2.24,-16.88) rectangle (7.5,2.24);
\draw [line width=1pt] (0,0)-- (5,0);
\draw [line width=1pt] (0,0)-- (0,-3);
\draw [line width=1pt] (0,-3)-- (1,-3);
\draw [line width=1pt] (1,-3)-- (1,-2);
\draw [line width=1pt] (1,-2)-- (3,-2);
\draw [line width=1pt] (3,-2)-- (3,-1);
\draw [line width=1pt] (3,-1)-- (5,-1);
\draw [line width=1pt] (5,-1)-- (5,0);
\draw [line width=1pt] (0,-2)-- (1,-2);
\draw [line width=1pt] (0,-1)-- (3,-1);
\draw [line width=1pt] (1,-2)-- (1,0);
\draw [line width=1pt] (2,-2)-- (2,0);
\draw [line width=1pt] (3,-1)-- (3,0);
\draw [line width=1pt] (4,-1)-- (4,0);
\draw [line width=1pt] (0,-2)-- (3,-5);
\draw [line width=1pt] (0,-1)-- (3.18,-4.16);
\draw [line width=1pt] (0,0)-- (3.56,-3.58);
\draw [line width=1pt] (1,0)-- (4,-3);
\draw [line width=1pt] (2,0)-- (4.58,-2.56);
\draw [line width=1pt] (3,0)-- (5,-2);
\draw [line width=1pt] (4,0)-- (6,-2);
\draw (3.58,-3.2) node[anchor=north west] {\small{$0$}};
\draw (3.18,-3.76) node[anchor=north west] {\small{$1$}};
\draw (3,-4.62) node[anchor=north west] {\small{$2$}};
\draw (3.74,-2.6) node[anchor=north west] {\small{$-1$}};
\draw (4.34,-2.2) node[anchor=north west] {\small{$-2$}};
\draw (4.78,-1.6) node[anchor=north west] {\small{$-3$}};
\draw (5.74,-1.6) node[anchor=north west] {\small{$-4$}};
\end{tikzpicture}
\vspace{-8cm}
\caption{Diagonals of the partition $\pi = (5,3,1)$}
\label{fig2}
\end{figure}
Note that the successive ranks given by Atkin \cite{A66} are the hook differences on the diagonal $0$.\quad\\\par Let $L,m,n$ be non-negative integers. We now recall some notations from the theory of $q$-series that can be found in \cite{A98}. To begin with, the conventional $q$-Pochhammer symbol is defined as
\begin{align*}
    (a)_L = (a;q)_L &:= \prod_{k=0}^{L-1}(1-aq^k),\\
    (a)_{\infty} = (a;q)_{\infty} &:= \lim_{L\rightarrow \infty}(a)_L\,\,\text{where}\,\,\lvert q\rvert<1.
\end{align*}\\
We define the $q$-binomial (Gaussian) coefficient as\\
\begin{align*}
    \left[\begin{matrix}m\\n\end{matrix}\right]_q := \Bigg\{\begin{array}{lr}
        \dfrac{(q)_m}{(q)_n(q)_{m-n}}\quad\text{for } m\ge n\ge 0,\\
        0\qquad\qquad\quad\text{otherwise}.\end{array}
\end{align*}\\
For $m, n\ge 0$, $\left[\begin{matrix}m+n\\n\end{matrix}\right]_q$ is the generating function for partitions into at most $n$ parts each of size at most $m$ (see \cite[Chapter $3$]{A98}).\quad\\\par Throughout the remainder of the paper, $P(q)\ge 0$ means that a power series in $q$, $P(q)$, has non-negative coefficients.\quad\\\par In \cite[Theorem $1$]{ABBBFV87} Andrews \textit{et al.} proved the following theorem using recurrences.\quad\\\par
\begin{theorem}\label{thm11}
The generating function $D_{K,i}(N,M;\alpha,\beta)$ of partitions with at most $M$ parts, largest part not exceeding $N$, and hook differences on the $(1-\beta)$th diagonal at least $\beta-i+1$ and on the $(\alpha-1)$th diagonal at most $K-\alpha-i-1$ is given by\\
\begin{equation}\label{eq11}
\begin{multlined}
D_{K,i}(N,M;\alpha,\beta;q) = D_{K,i}(N,M;\alpha,\beta)\\ = \sum\limits_{j\in\mathbb{Z}}\Bigg\{q^{j((\alpha+\beta)Kj+K\beta-(\alpha+\beta)i)}\left[\begin{matrix}M+N\\M-Kj\end{matrix}\right]_{q}-q^{((\alpha+\beta)j+\beta)(Kj+i)}\left[\begin{matrix}M+N\\M-Kj-i\end{matrix}\right]_{q}\Bigg\}.\\
\end{multlined}
\end{equation}
Here the following conditions apply: $\alpha,\beta\in\mathbb{N}\cup\{0\}$, $0 < i < K$, and $\beta-i\le N-M\le K-\alpha-i$ with the added restrictions that the largest part exceeds $M-i$ if $\beta = 0$ and the number of parts exceeds $N+i$ if $\alpha = 0$.\\    
\end{theorem}
\begin{remark}\label{rmk1}
In \cite[Theorem $1$]{ABBBFV87}, $i$ was confined to the range $1\le i\le K/2$. However, it can be extended to the larger range $1\le i\le K-1$ as stated in Theorem \ref{thm11} by using simple conjugation argument.\\    
\end{remark}
Since $D_{K,i}(N,M;\alpha,\beta)$ is a generating function of partitions, we have an immediate Corollary to Theorem \ref{thm11} which is as follows.\\
\begin{cor}\label{cor12}
Let $K,i$ be positive integers such that $0 < i < K$ and $N,M,\alpha,\beta$ be non-negative integers such that $1\le \alpha+\beta\le K-1$ and $\beta-i\le N-M\le K-\alpha-i$. Then, $D_{K,i}(N,M;\alpha,\beta)\ge 0.$\\
\end{cor}
\begin{remark}\label{rmk2}
Note that, one can verify the following symmetry\\
\begin{align*}
D_{K,i}(N,M;\alpha,\beta) = D_{K,K-i}(M,N;\beta,\alpha).\numberthis\label{eq12}\\
\end{align*}
\end{remark}
\begin{remark}\label{rmk3}
The symmetry in \eqref{eq12} holds for $D_{K,i}(N,M;\alpha,\beta)$ too where $\alpha$ and $\beta$ are non-integral but $\alpha K$, $\beta K$, $\alpha i$, and $\beta i$ are integers.
\end{remark}\par\quad\\
Define $G(N,M;\alpha,\beta,K) = D_{2K,K}(N,M;\alpha,\beta)$. Then, we have \cite[eq. ($1.5$)]{W01},\\ 
\begin{align*}
G(N,M;\alpha,\beta,K) = \sum\limits_{j\in\mathbb{Z}}(-1)^jq^{\frac{1}{2}Kj((\alpha+\beta)j+\alpha-\beta)}\left[\begin{matrix}M+N\\N-Kj\end{matrix}\right]_{q}.\numberthis\label{eq13}\\    
\end{align*}
Now, note that one of the mod $3$ conjectures due to Borwein \cite{A95} can be stated as\\
\begin{equation*}
G(N,N;4/3,5/3,3)\ge 0,
\end{equation*}
\begin{equation*}
G(N+1,N-1;2/3,7/3,3)\ge 0,
\end{equation*}
and
\begin{equation*}
G(N+1,N-1;1/3,8/3,3)\ge 0.
\end{equation*}
\\\par All the three inequalities above were proven recently by Wang \cite{W22}.\\\par In \cite[Conjecture $6$]{B96}, Bressoud made the following interesting generalization.\\
\begin{conj}\label{conj13}
Let $K$ be a positive integer and $N,M,\alpha K,\beta K$ be non-negative integers such that $1\le \alpha+\beta\le 2K-1$ (strict inequalities when $K = 2$) and $\beta-K\le N-M\le K-\alpha$. Then, $G(N,M;\alpha,\beta,K)$ is a polynomial in $q$ with non-negative coefficients.\\
\end{conj}
Many cases of Conjecture \ref{conj13} were proven in the literature \cite{B20,B22,BW05,B81,IKS99,W22,W01,W03}.\\

The rest of the paper is organized as follows. In Section \ref{s2}, we present the statement of the main results where we state a new generalization of Conjecture \ref{conj13}, namely, Conjecture \ref{conj21} and propose new infinite hierarchies of families of non-negative $q$-series identities such as in Theorem \ref{thm26} and Theorem \ref{thm27}. In Section \ref{s3}, we state five elegant and important positivity-preserving transformations for $q$-binomial coefficients. In Section \ref{s4}, we present the proofs of the non-negativity of the infinite families stated in Section \ref{s2}.

\section{Main Results}\label{s2}
In this section, we first state a new extension of Conjecture \ref{conj13} which is much more general and is as follows.\\
\begin{conj}\label{conj21}
Let $K,i$ be positive integers such that $0 < i < K$ and $N,M,\alpha K,\beta K,\alpha i,\beta i$ be non-negative integers such that $1\le \alpha+\beta\le K-1$ (strict inequalities when $K = 4$ and $i=2$) and $\beta-i\le N-M\le K-\alpha-i$. Then, $D_{K,i}(N,M;\alpha,\beta)$, as defined in Theorem \ref{thm11}, is a polynomial in $q$ with non-negative coefficients.\\
\end{conj}
It can be easily verified that Conjecture \ref{conj13} is the special case $(i,K) = (K,2K)$ of Conjecture \ref{conj21}.\\\par Now, we state two examples of polynomials having non-negative coefficients for which all the conditions stated above in Conjecture \ref{conj21} are satisfied. Let $a\in\{0,1\}$. Then, for any non-negative integer $N$, we have\\
\begin{enumerate}[(i)]
\item \begin{align*}
D_{6,2}(N+a,N;1,1/2) &:= {\displaystyle\sum\limits_{j\in\mathbb{Z}}\left(\frac{j+1}{3}\right)q^{j^2}\left[\begin{matrix}2N+a\\N-2j\end{matrix}\right]_{q}}\\ &= D_{6,4}(N,N+a;1/2,1)\ge 0,    
\end{align*}
\item \begin{align*}
D_{9,3}(N,N+a;2/3,1/3) &:= {\displaystyle\sum\limits_{j\in\mathbb{Z}}\left(\frac{j+1}{3}\right)q^{j^2}\left[\begin{matrix}2N+a\\N+3j\end{matrix}\right]_{q}}\\ &= {\displaystyle\sum\limits_{j\in\mathbb{Z}}\left(\frac{-j+1}{3}\right)q^{j^2}\left[\begin{matrix}2N+a\\N-3j\end{matrix}\right]_{q}}\\ &= D_{9,6}(N+a,N;1/3,2/3)\ge 0,\\    
\end{align*}
\end{enumerate}
where the last equality in both (i) and (ii) follows from the symmetry relation \eqref{eq12} and we define the usual Legendre symbol (mod $3$) as\\ $$\left(\frac{j}{3}\right) = \begin{cases}
1,\quad\quad\text{if}\,\,\,j\equiv 1 \pmod{3},\\
-1,\quad\,\text{if}\,\,\,j\equiv 2 \pmod{3},\\
0,\quad\quad\text{if}\,\,\,j\equiv 0 \pmod{3}.
\end{cases}$$\\\par We verified (i) and (ii) using extensive computer checks.\\
\begin{theorem}\label{thm22}
For any non-negative integer $L$ and positive integers $p^{\prime}$, $p$ such that $p^{\prime} > p$,\\
\begin{equation}\label{eq21}
\begin{multlined}
D_{p^{\prime},s}\Bigg(\left\lceil\frac{L+r-s}{2}\right\rceil,\left\lfloor\frac{L-r+s}{2}\right\rfloor;p-r,r\Bigg)\\ = \sum\limits_{j\in\mathbb{Z}}\Bigg\{q^{j^2pp^{\prime}+(rp^{\prime}-sp)j}\left[\begin{matrix}L\\\left\lfloor\frac{L-r+s}{2}\right\rfloor-jp^{\prime}\end{matrix}\right]_{q}-q^{(jp+r)(jp^{\prime}+s)}\left[\begin{matrix}L\\\left\lfloor\frac{L-r-s}{2}\right\rfloor-jp^{\prime}\end{matrix}\right]_{q}\Bigg\}\ge 0,\\   
\end{multlined}    
\end{equation}
where $r$ and $s$ are integers such that $0 < r < p$ and $0 < s < p^{\prime}$.\\
\end{theorem}
\par Observe that Theorem $\ref{thm22}$ follows immediately from Corollary \ref{cor12} and is the starting point of proofs of all positivity results below.\\
\begin{theorem}\label{thm23}
For all non-negative integers $L$, positive integers $p,p^{\prime}$ such that $p^{\prime} > p$, and any integer $n\ge 2$,\\
\begin{equation}\label{eq22}
\begin{multlined}
D_{2^np^{\prime},2^ns}\Bigg(L+2^{n-2}+2^{n-1}(r-s),L-2^{n-2}-2^{n-1}(r-s);\\\frac{\frac{2^{2n}-1}{3}(2p^{\prime}-1)-\frac{2^{2n+1}+1}{3}r+p}{2^n},\frac{\frac{2^{2n+1}+1}{3}r+\frac{2^{2n}-1}{3}}{2^n}\Bigg) \ge 0,\\   
\end{multlined}    
\end{equation}
where $r$, $s$ are integers such that $0 < r < p$ and $0 < s < p^{\prime}$.\\
\end{theorem}
\par Note that Theorem \ref{thm23} is consistent with Conjecture \ref{conj21}. We now have an immediate Corollary of Theorem \ref{thm23} which is as follows.\\
\begin{cor}\label{cor24}
Note that, for even integers $p^{\prime}$ and $s = p^{\prime}/2$ in Theorem \ref{thm23}, we get special cases of Bressoud's Conjecture \ref{conj13}. For instance, for $p^{\prime} = 2\Tilde{p}$ and $s = \Tilde{p}$ and for all integers $n\ge 2$, we have\\
\begin{equation}\label{eq23}
\begin{multlined}
G\Bigg(L+2^{n-2}+2^{n-1}(r-\Tilde{p}),L-2^{n-2}-2^{n-1}(r-\Tilde{p});\\\frac{\frac{2^{2n}-1}{3}(4\Tilde{p}-1)-\frac{2^{2n+1}+1}{3}r+p}{2^n},\frac{\frac{2^{2n+1}+1}{3}r+\frac{2^{2n}-1}{3}}{2^n},2^n\Tilde{p}\Bigg)\ge 0.    
\end{multlined}
\end{equation}
where $0 < p < 2\Tilde{p}$ and $0 < r < p$.\\
\end{cor}
\begin{theorem}\label{thm25}
For $L\in\mathbb{N}$, $\nu\in\mathbb{Z}_{> 0}$, $s = 0,1,\ldots,\nu-1$, and all integers $n\ge 2$,
\begin{equation}\label{eq24}
\begin{multlined}
G\Bigg(L-2^{n-1}s-2^{n-2},L+2^{n-1}s+2^{n-2};\\\frac{(\frac{2^{2n}-1}{3}\nu+\frac{2^{2n-1}+1}{3})(\nu+s+1)}{2^{n-2}(2\nu+1)},\frac{(\frac{2^{2n}-1}{3}\nu+\frac{2^{2n-1}+1}{3})(\nu-s)}{2^{n-2}(2\nu+1)},2^{n-1}(2\nu+1)\Bigg)\ge 0.
\end{multlined}
\end{equation}
\end{theorem}
\par Note that Theorem \ref{thm25} is consistent with Conjecture \ref{conj13}.\\
\par Next, we state a Corollary of Theorem \ref{thm23} which is as follows.\\
\begin{theorem}\label{thm26}
For all non-negative integers $L,t$, positive integers $p,p^{\prime}$ such that $p^{\prime} > p$, and any integer $n\ge 2$,
\begin{equation}\label{eq25}
\begin{multlined}
D_{3^t\cdot2^np^{\prime},3^t\cdot2^ns}\Bigg(L+3^t\cdot2^{n-2}+3^t\cdot2^{n-1}(r-s),L-3^t\cdot2^{n-2}-3^t\cdot2^{n-1}(r-s);\\\frac{\frac{2^{2n}-1}{3}(2p^{\prime}-1)-\frac{2^{2n+1}+1}{3}r+p}{2^n}-(3^t-1)\cdot2^{n-2}+(3^t-1)\cdot2^{n-1}(p^{\prime}-r),\\\frac{\frac{2^{2n+1}+1}{3}r+\frac{2^{2n}-1}{3}}{2^n}+(3^t-1)\cdot2^{n-2}+(3^t-1)\cdot2^{n-1}r\Bigg) \ge 0,   
\end{multlined}    
\end{equation}
where $r$, $s$ are integers such that $0 < r < p$ and $0 < s < p^{\prime}$.\\
\end{theorem}
\par Observe that, Theorem \ref{thm26} is consistent with Conjecture \ref{conj21} and for $t = 0$, Theorem \ref{thm26} becomes Theorem \ref{thm23}. Finally, we state another result which is consistent with Conjecture \ref{conj21}.\\
\begin{theorem}\label{thm27}
For all non-negative integers $L,t$, positive integers $p,p^{\prime}$ such that $p^{\prime} > p$, and any integer $n\ge 2$,
\begin{equation}\label{eq26}
\begin{multlined}
D_{3^t\cdot2^np^{\prime},3^t\cdot2^ns}\Bigg(L+3^t\cdot2^{n-2}+3^t\cdot2^{n-1}(r-s),L-3^t\cdot2^{n-2}-3^t\cdot2^{n-1}(r-s);\\\frac{\left(\frac{2^{2n+1}-8}{3}\cdot3^t+4t+2\right)p^{\prime}-\frac{2^{2n}-1}{3}\cdot3^t-\left(\frac{2^{2n+1}-2}{3}\cdot3^t+1\right)r+p}{2^n},\\\frac{\left(\frac{2^{2n+1}-2}{3}\cdot3^t+1\right)r+\frac{2^{2n}-1}{3}\cdot3^t}{2^n}\Bigg) \ge 0,   
\end{multlined}    
\end{equation}
where $r$, $s$ are integers such that $0 < r < p$ and $0 < s < p^{\prime}$.\\
\end{theorem}

\section{Some Positivity-Preserving Transformations}\label{s3}
In this section, we state some positivity-preserving transformations for certain $q$-binomial coefficients.\\
\begin{theorem}(Berkovich \cite[Theorem $2.1$]{B20,B22})\label{thm31}
For $L\in\mathbb{N}$ and $a\in\mathbb{Z}$, we have\\
\begin{equation}\label{eq31}
\sum\limits_{k\ge 0}C_{L,k}(q)\left[\begin{matrix}k \\ \left\lfloor\frac{k-a}{2}\right\rfloor\end{matrix}\right]_{q} = q^{T(a)}\left[\begin{matrix}2L+1 \\ L-a\end{matrix}\right]_{q},    
\end{equation}\\
where $T(j) := \binom{j+1}{2}$ and\\
\begin{equation}\label{eq32}
C_{L,k}(q) = \sum\limits_{m=0}^{L}q^{T(m)+T(m+k)}\left[\begin{matrix}L \\ m,k\end{matrix}\right]_{q},\\
\end{equation}
where\\
\begin{equation}\label{eq33}
\left[\begin{matrix}L \\ m,k\end{matrix}\right]_{q} = \left[\begin{matrix}L \\ m\end{matrix}\right]_{q}\left[\begin{matrix}L-m \\ k\end{matrix}\right]_{q} = \left[\begin{matrix}L \\ k\end{matrix}\right]_{q}\left[\begin{matrix}L-k \\ m\end{matrix}\right]_{q}\ge 0.
\end{equation}\\
\end{theorem}
Now, note that $C_{L,k}(q)\ge 0$. It is then easy to verify that for any identity of the form\\
\begin{equation}\label{eq34}
F_C(L,q) = \sum\limits_{j\in\mathbb{Z}}\alpha(j,q)\left[\begin{matrix}L \\ \left\lfloor\frac{L-j}{2}\right\rfloor\end{matrix}\right]_{q},    
\end{equation}\\
using transformation \eqref{eq31}, the following identity holds\\
\begin{equation}\label{eq35}
\sum\limits_{k\ge 0}C_{L,k}(q)F_C(k,q) = \sum\limits_{j\in\mathbb{Z}}\alpha(j,q)\sum\limits_{k\ge 0}C_{L,k}(q)\left[\begin{matrix}k \\ \left\lfloor\frac{k-j}{2}\right\rfloor\end{matrix}\right]_{q} = \sum\limits_{j\in\mathbb{Z}}\alpha(j,q)q^{T(j)}\left[\begin{matrix}2L+1 \\ L-j\end{matrix}\right]_{q}.    
\end{equation}\\
Hence, if $F_C(L,q)\ge 0$, then\\
\begin{equation}\label{eq36}
\sum\limits_{j\in\mathbb{Z}}\alpha(j,q)q^{T(j)}\left[\begin{matrix}2L+1 \\ L-j\end{matrix}\right]_{q}\ge 0.    
\end{equation}\\
So, we say that the transformation \eqref{eq31} is positivity-preserving.\\
\begin{theorem}(Berkovich-Warnaar \cite{BW05})\label{thm32}
For $L\in\mathbb{N}$ and $a\in\mathbb{Z}$, we have\\
\begin{equation}\label{eq37}
\sum\limits_{k\ge 0}O_{L,k}(q)\left[\begin{matrix}2k+1 \\ k-a\end{matrix}\right]_{q} = q^{4T(a)}\left[\begin{matrix}2L \\ L-2a-1\end{matrix}\right]_{q},    
\end{equation}\\
where\\
\begin{equation}\label{eq38}
O_{L,k}(q) = \sum\limits_{m=0}^{L}q^{2T(k)+2T(m+k)}\left[\begin{matrix}L \\ m,2k+1\end{matrix}\right]_{q}\ge 0.
\end{equation}\\
\end{theorem}
It is then easy to verify that for any identity of the form\\
\begin{equation}\label{eq39}
F_O(L,q) = \sum\limits_{j\in\mathbb{Z}}\alpha(j,q)\left[\begin{matrix}2L+1 \\ L-j\end{matrix}\right]_{q},    
\end{equation}\\
using transformation \eqref{eq37}, the following identity holds\\
\begin{equation}\label{eq310}
\sum\limits_{k\ge 0}O_{L,k}(q)F_O(k,q) = \sum\limits_{j\in\mathbb{Z}}\alpha(j,q)\sum\limits_{k\ge 0}O_{L,k}(q)\left[\begin{matrix}2k+1 \\ k-j\end{matrix}\right]_{q} = \sum\limits_{j\in\mathbb{Z}}\alpha(j,q)q^{4T(j)}\left[\begin{matrix}2L \\ L-2j-1\end{matrix}\right]_{q}.    
\end{equation}\\
Hence, if $F_O(L,q)\ge 0$, then\\
\begin{equation}\label{eq311}
\sum\limits_{j\in\mathbb{Z}}\alpha(j,q)q^{4T(j)}\left[\begin{matrix}2L \\ L-2j-1\end{matrix}\right]_{q}\ge 0.    
\end{equation}\\
So, again transformation \eqref{eq37} is positivity-preserving.\\
\begin{theorem}(Warnaar \cite[Corollary $2.6$]{W03})\label{thm33}
For $L\in\mathbb{N}$ and $a\in\mathbb{Z}$, we have\\
\begin{equation}\label{eq312}
\sum\limits_{k\ge 0}W_{L,k}(q)\left[\begin{matrix}2k \\ k-a\end{matrix}\right]_{q} = q^{2a^2}\left[\begin{matrix}2L \\ L-2a\end{matrix}\right]_{q},    
\end{equation}\\
where\\
\begin{equation}\label{eq313}
W_{L,k}(q) = \sum\limits_{m=0}^{L}q^{k^2+(m+k)^2}\left[\begin{matrix}L \\ m,2k\end{matrix}\right]_{q}\ge 0.
\end{equation}\\
\end{theorem}
It is then easy to verify that for any identity of the form\\
\begin{equation}\label{eq314}
F_W(L,q) = \sum\limits_{j\in\mathbb{Z}}\alpha(j,q)\left[\begin{matrix}2L \\ L-j\end{matrix}\right]_{q},    
\end{equation}\\
using transformation \eqref{eq312}, the following identity holds\\
\begin{equation}\label{eq315}
\sum\limits_{k\ge 0}W_{L,k}(q)F_W(k,q) = \sum\limits_{j\in\mathbb{Z}}\alpha(j,q)\sum\limits_{k\ge 0}W_{L,k}(q)\left[\begin{matrix}2k \\ k-j\end{matrix}\right]_{q} = \sum\limits_{j\in\mathbb{Z}}\alpha(j,q)q^{2j^2}\left[\begin{matrix}2L \\ L-2j\end{matrix}\right]_{q}.    
\end{equation}\\
Hence, if $F_W(L,q)\ge 0$, then\\
\begin{equation}\label{eq316}
\sum\limits_{j\in\mathbb{Z}}\alpha(j,q)q^{2j^2}\left[\begin{matrix}2L \\ L-2j\end{matrix}\right]_{q}\ge 0.    
\end{equation}\\
So, again transformation \eqref{eq312} is positivity-preserving.\\
\begin{theorem}(\cite[Lemma $2.6$]{BW05} $L$, $j$, $r$ even)\label{thm34}
For integers $L$ and $j$, we have\\
\begin{equation}\label{eq317}
\sum\limits_{r=0}^{\left\lfloor\frac{L}{3}\right\rfloor}A_{L,r}(q)\left[\begin{matrix}2r \\ r-j\end{matrix}\right]_{q^3} = q^{3j^2}\left[\begin{matrix}2L \\ L-3j\end{matrix}\right]_{q},
\end{equation}\\
where\\
\begin{equation}\label{eq318}
A_{L,r}(q) = \dfrac{q^{3r^2}(q^3;q^3)_{L-r-1}(1-q^{2L})}{(q^3;q^3)_{2r}(q;q)_{L-3r}}.    
\end{equation}\\
\end{theorem}
Berkovich and Warnaar \cite{BW05} showed that\\
\begin{equation}\label{eq319}
f_{L,r}(q) = \dfrac{(q^3;q^3)_{\frac{1}{2}(L-r-2)}(1-q^L)}{(q^3;q^3)_r(q;q)_{\frac{1}{2}(L-3r)}}    
\end{equation}\\ 
is a  polynomial with non-negative coefficients for $0\le 3r\le L$ and $r\equiv L\pmod 2$. It is then evident from \eqref{eq318} and \eqref{eq319} that\\ $$A_{L,r}(q) = q^{3r^2}f_{2L,2r}(q)\\$$ has non-negative coefficients.\par\quad\par
It is then easy to verify that for any identity of the form\\
\begin{equation}\label{eq320}
F_A(L,q) = \sum\limits_{j\in\mathbb{Z}}\alpha(j,q)\left[\begin{matrix}2L \\ L-j\end{matrix}\right]_{q^3},    
\end{equation}\\
using transformation \eqref{eq317}, the following identity holds\\
\begin{equation}\label{eq321}
\sum\limits_{r\ge 0}A_{L,r}(q)F_A(r,q) = \sum\limits_{j\in\mathbb{Z}}\alpha(j,q)\sum\limits_{r\ge 0}A_{L,r}(q)\left[\begin{matrix}2r \\ r-j\end{matrix}\right]_{q^3} = \sum\limits_{j\in\mathbb{Z}}\alpha(j,q)q^{3j^2}\left[\begin{matrix}2L \\ L-3j\end{matrix}\right]_{q}.    
\end{equation}\\
Hence, if $F_A(L,q)\ge 0$, then\\
\begin{equation}\label{eq322}
\sum\limits_{j\in\mathbb{Z}}\alpha(j,q)q^{3j^2}\left[\begin{matrix}2L \\ L-3j\end{matrix}\right]_{q}\ge 0.    
\end{equation}\\
So, again transformation \eqref{eq317} is positivity-preserving.\\
\begin{theorem}(\cite[Lemma $2.6$]{BW05} $L$, $j$, $r$ odd)\label{thm35}
For integers $L$ and $j$, we have\\
\begin{equation}\label{eq323}
\sum\limits_{r=0}^{\left\lfloor\frac{L}{3}\right\rfloor}\Tilde{A}_{L,r}(q)\left[\begin{matrix}2r+1 \\ r-j\end{matrix}\right]_{q^3} = q^{3j^2+3j}\left[\begin{matrix}2L+1 \\ L-3j-1\end{matrix}\right]_{q},
\end{equation}\\
where\\
\begin{equation}\label{eq324}
\Tilde{A}_{L,r}(q) = \dfrac{q^{3r^2+3r}(q^3;q^3)_{L-r-1}(1-q^{2L+1})}{(q^3;q^3)_{2r+1}(q;q)_{L-3r-1}}.    
\end{equation}\\
\end{theorem}
It is then evident from \eqref{eq319} and \eqref{eq324} that\\ $$\Tilde{A}_{L,r}(q) = q^{3r^2+3r}f_{2L+1,2r+1}(q)\\$$ has non-negative coefficients.\par\quad\par
It is then easy to verify that for any identity of the form\\
\begin{equation}\label{eq325}
F_{\Tilde{A}}(L,q) = \sum\limits_{j\in\mathbb{Z}}\alpha(j,q)\left[\begin{matrix}2L+1 \\ L-j\end{matrix}\right]_{q^3},    
\end{equation}\\
using transformation \eqref{eq323}, the following identity holds\\
\begin{equation}\label{eq326}
\sum\limits_{r\ge 0}\Tilde{A}_{L,r}(q)F_{\Tilde{A}}(r,q) = \sum\limits_{j\in\mathbb{Z}}\alpha(j,q)\sum\limits_{r\ge 0}\Tilde{A}_{L,r}(q)\left[\begin{matrix}2r+1 \\ r-j\end{matrix}\right]_{q^3} = \sum\limits_{j\in\mathbb{Z}}\alpha(j,q)q^{3j^2+3j}\left[\begin{matrix}2L+1 \\ L-3j-1\end{matrix}\right]_{q}.    
\end{equation}\\
Hence, if $F_{\Tilde{A}}(L,q)\ge 0$, then\\
\begin{equation}\label{eq327}
\sum\limits_{j\in\mathbb{Z}}\alpha(j,q)q^{3j^2+3j}\left[\begin{matrix}2L+1 \\ L-3j-1\end{matrix}\right]_{q}\ge 0.    
\end{equation}\\
So, again transformation \eqref{eq323} is positivity-preserving.\\

\section{Proofs of some Main Results}\label{s4}
In this section, we prove Theorem \ref{thm23}, Theorem \ref{thm25}, Theorem \ref{thm26}, and Theorem \ref{thm27}.\\
\subsection{Proof of Theorem \ref{thm23}}\label{ss41}
Applying \eqref{eq35} to \eqref{eq21}, we get\\
\begin{equation}\label{eq41}
D_{2p^{\prime},2s}\left(L+1+r-s, L-r+s; \frac{2p^{\prime}-1-3r+p}{2}, \frac{3r+1}{2}\right)\ge 0.    
\end{equation}\\
Now, applying \eqref{eq310} to \eqref{eq41}, we get\\
\begin{equation}\label{eq42}
D_{4p^{\prime},4s}\left(L+1+2r-2s, L-1-2r+2s; \frac{10p^{\prime}-5-11r+p}{4}, \frac{11r+5}{4}\right)\ge 0.    
\end{equation}\\
Now, consider
\begin{equation}\label{eq43}
\begin{multlined}
D_{K,i}\left(L+a,L-a;\alpha,\beta\right) = \sum\limits_{j\in\mathbb{Z}}\Bigg\{q^{j(K(\alpha+\beta)j+K\beta-(\alpha+\beta)i)}\left[\begin{matrix}2L \\ L-a-Kj\end{matrix}\right]_{q}\\ - q^{((\alpha+\beta)j+\beta)(Kj+i)}\left[\begin{matrix}2L \\ L-a-Kj-i\end{matrix}\right]_{q}\Bigg\}.    
\end{multlined}
\end{equation}
Applying \eqref{eq315} to \eqref{eq43}, we get\\
\begin{equation}\label{eq44}
\begin{multlined}
\sum\limits_{r\ge 0}W_{L,r}(q)D_{K,i}(r+a,r-a;\alpha,\beta)\\ = q^{2a^2}D_{2K,2i}\left(L+2a,L-2a;\dfrac{\alpha-4a+2K-2i}{2},\dfrac{\beta+4a+2i}{2}\right).\\    
\end{multlined}
\end{equation}
Now, recall that $W_{L,r}(q)\ge 0$ and so if\\ 
\begin{equation}\label{eq45}
D_{K,i}(L+a,L-a;\alpha,\beta)\ge 0,    
\end{equation}\\ 
then\\
\begin{equation}\label{eq46}
D_{2K,2i}\left(L+2a,L-2a;\dfrac{\alpha-4a+2K-2i}{2},\dfrac{\beta+4a+2i}{2}\right)\ge 0.    
\end{equation}\\
Now, iterating the same procedure $(n-2)$ times for all $n\ge 2$, we have\\
\begin{equation}\label{eq47}
\begin{multlined}
D_{2^{n-2}K,2^{n-2}i}\Bigg(L+2^{n-2}a,L-2^{n-2}a;\\\dfrac{\alpha-\dfrac{4}{3}(4^{n-2}-1)a+\dfrac{2}{3}(4^{n-2}-1)(K-i)}{2^{n-2}},\dfrac{\beta+\dfrac{4}{3}(4^{n-2}-1)a+\dfrac{2}{3}(4^{n-2}-1)i}{2^{n-2}}\Bigg)\ge 0\\
\end{multlined}    
\end{equation}
if \eqref{eq45} holds. In particular, together with \eqref{eq42}, \eqref{eq47} implies Theorem \ref{thm23}.\qed\\
\subsection{Proof of Theorem \ref{thm25}}\label{ss42}
Let\\
\begin{equation*}
F(L,\nu,s,q) = \sum\limits_{j\in\mathbb{Z}}(-1)^jq^{(\nu+1)(2\nu+1)j^2+(\nu+1)(2s+1)j}\left[\begin{matrix}2L+1 \\ L-s-(2\nu+1)j\end{matrix}\right]_{q}.    
\end{equation*}\\ 
Clearly, observe that\\
\begin{equation}\label{eq48}
F(L,\nu,s,q) = G\Bigg(L-s,L+s+1;\\\frac{2(\nu+1)(\nu+s+1)}{2\nu+1},\frac{2(\nu+1)(\nu-s)}{2\nu+1},2\nu+1\Bigg)\ge 0.   \end{equation}\\
Also, note that $q^{T(s)}F(L,\nu,s,q)$ is the right-hand side of \cite[eq. ($2.23$)]{B20} and so, $F(L,\nu,s,q)\ge 0$. Now, applying \eqref{eq310} to \eqref{eq48}, we get\\
\begin{equation}\label{eq49}
G\Bigg(L-2s-1,L+2s+1;\\\frac{(5\nu+3)(\nu+s+1)}{2\nu+1},\frac{(5\nu+3)(\nu-s)}{2\nu+1},4\nu+2\Bigg)\ge 0.    
\end{equation}\\
Now, applying \eqref{eq315} to \eqref{eq49} and iterating using \eqref{eq315} $(n-2)$ times for all $n\ge 2$, we get Theorem \ref{thm25} and note that \eqref{eq49} is the $n=2$ case of \eqref{eq24}.\qed\\
\subsection{Proof of Theorem \ref{thm26}}\label{ss43}
Consider
\begin{equation}\label{eq410}
\begin{multlined}
D_{K,i}\left(L+a,L-a;\alpha,\beta,q\right) = \sum\limits_{j\in\mathbb{Z}}\Bigg\{q^{j(K(\alpha+\beta)j+K\beta-(\alpha+\beta)i)}\left[\begin{matrix}2L \\ L-a-Kj\end{matrix}\right]_{q}\\ - q^{((\alpha+\beta)j+\beta)(Kj+i)}\left[\begin{matrix}2L \\ L-a-Kj-i\end{matrix}\right]_{q}\Bigg\}.    
\end{multlined}
\end{equation}
Making the substitution $q\longrightarrow q^3$ in \eqref{eq410} and then applying \eqref{eq321}, we get\\
\begin{equation}\label{eq411}
\begin{multlined}
\sum\limits_{r=0}^{\left\lfloor\frac{L}{3}\right\rfloor}A_{L,r}(q)D_{K,i}(r+a,r-a;\alpha,\beta,q^3)\\ = q^{3a^2}D_{3K,3i}(L+3a,L-3a;\alpha-2a+K-i,\beta+2a+i,q).\\    
\end{multlined}
\end{equation}
Now, recall that $A_{L,r}(q)\ge 0$ and so if\\ 
\begin{equation}\label{eq412}
D_{K,i}(L+a,L-a;\alpha,\beta,q)\ge 0,    
\end{equation}\\ 
then\\
\begin{equation}\label{eq413}
D_{3K,3i}(L+3a,L-3a;\alpha-2a+K-i,\beta+2a+i,q)\ge 0.    
\end{equation}\\
Now, iterating the same process $t$ times for all $t\ge 0$, we have\\
\begin{equation}\label{eq414}
\begin{multlined}
D_{3^tK,3^ti}\left(L+3^ta,L-3^ta;\alpha-(3^t-1)a+\dfrac{3^t-1}{2}(K-i),\beta+(3^t-1)a+\dfrac{3^t-1}{2}i,q\right)\ge 0\\
\end{multlined}    
\end{equation}
if \eqref{eq412} holds. In particular, together with Theorem \ref{thm23}, \eqref{eq414} implies Theorem \ref{thm26}.\qed\\
\subsection{Proof of Theorem \ref{thm27}}\label{ss44}
Making the substitution $q\longrightarrow q^3$ in \eqref{eq41} and then applying \eqref{eq326}, we get\\
\begin{equation}\label{eq415}
\begin{multlined}
\sum\limits_{\Tilde{r}=0}^{\left\lfloor\frac{L}{3}\right\rfloor}\Tilde{A}_{L,\Tilde{r}}(q)D_{2p^{\prime},2s}\left(\Tilde{r}+1+r-s,\Tilde{r}-r+s;\frac{2p^{\prime}-1-3r+p}{2},\frac{3r+1}{2},q^3\right)\\ = q^{3(r-s)^2+3(r-s)}D_{6p^{\prime},6s}\left(L+2+3r-3s,L-1-3r+3s;\frac{6p^{\prime}-3-7r+p}{2},\frac{7r+3}{2},q\right).\\    
\end{multlined}
\end{equation}
Now, recall that $\Tilde{A}_{L,\Tilde{r}}(q)\ge 0$ and thus, by \eqref{eq41}, we have\\
\begin{equation}\label{eq416}
D_{6p^{\prime},6s}\left(L+2+3r-3s,L-1-3r+3s;\frac{6p^{\prime}-3-7r+p}{2},\frac{7r+3}{2},q\right)\ge 0.    
\end{equation}\\
Now, iterating the same process $t$ times for all $t\ge 0$, we have\\
\begin{equation}\label{eq417}
\begin{multlined}
D_{3^t\cdot2p^{\prime},3^t\cdot2s}\Bigg(L+\frac{3^t+1}{2}+3^t(r-s),L-\frac{3^t-1}{2}-3^t(r-s);\\\frac{(4t+2)p^{\prime}-3^t-(2\cdot3^t+1)r+p}{2},\frac{(2\cdot3^t+1)r+3^t}{2}\Bigg)\ge 0.\\    
\end{multlined}    
\end{equation}
Next, applying \eqref{eq310} to \eqref{eq417}, we get\\
\begin{equation}\label{eq418}
\begin{multlined}
\sum\limits_{k\ge0}O_{L,k}(q)D_{3^t\cdot2p^{\prime},3^t\cdot2s}\Bigg(k+\frac{3^t+1}{2}+3^t(r-s),k-\frac{3^t-1}{2}-3^t(r-s);\\\frac{(4t+2)p^{\prime}-3^t-(2\cdot3^t+1)r+p}{2},\frac{(2\cdot3^t+1)r+3^t}{2}\Bigg)\\ = q^{2\left(\frac{3^t-1}{2}+3^t(r-s)\right)^2+2\left(\frac{3^t-1}{2}+3^t(r-s)\right)}D_{3^t\cdot4p^{\prime},3^t\cdot4s}\Bigg(L+3^t+3^t\cdot2(r-s),L-3^t-3^t\cdot2(r-s);\\\frac{(8\cdot3^t+4t+2)p^{\prime}-5\cdot3^t-(10\cdot3^t+1)r+p}{4},\frac{(10\cdot3^t+1)r+5\cdot3^t}{4}\Bigg).\\    
\end{multlined}
\end{equation}
Now, recall that $O_{L,k}(q)\ge 0$ and thus, by \eqref{eq417}, we have
\begin{equation}\label{eq419}
\begin{multlined}
D_{3^t\cdot4p^{\prime},3^t\cdot4s}\Bigg(L+3^t+3^t\cdot2(r-s),L-3^t-3^t\cdot2(r-s);\\\frac{(8\cdot3^t+4t+2)p^{\prime}-5\cdot3^t-(10\cdot3^t+1)r+p}{4},\frac{(10\cdot3^t+1)r+5\cdot3^t}{4}\Bigg)\ge 0.
\end{multlined}
\end{equation}
Next, applying \eqref{eq315} to \eqref{eq419}, we get
\begin{equation}\label{eq420}
\begin{multlined}
\sum\limits_{k\ge 0}W_{L,k}(q)D_{3^t\cdot4p^{\prime},3^t\cdot4s}\Bigg(k+3^t+3^t\cdot2(r-s),k-3^t-3^t\cdot2(r-s);\\\frac{(8\cdot3^t+4t+2)p^{\prime}-5\cdot3^t-(10\cdot3^t+1)r+p}{4},\frac{(10\cdot3^t+1)r+5\cdot3^t}{4}\Bigg)\\ = q^{2\left(3^t+3^t\cdot2(r-s)\right)^2}D_{3^t\cdot8p^{\prime},3^t\cdot8s}\Bigg(L+3^t\cdot2+3^t\cdot4(r-s),L-3^t\cdot2-3^t\cdot4(r-s);\\\frac{(40\cdot3^t+4t+2)p^{\prime}-21\cdot3^t-(42\cdot3^t+1)r+p}{8},\frac{(42\cdot3^t+1)r+21\cdot3^t}{8}\Bigg).   
\end{multlined}    
\end{equation}
Now, recall that $W_{L,k}(q)\ge 0$ and thus, by \eqref{eq419}, we have
\begin{equation}\label{eq421}
\begin{multlined}
D_{3^t\cdot8p^{\prime},3^t\cdot8s}\Bigg(L+3^t\cdot2+3^t\cdot4(r-s),L-3^t\cdot2-3^t\cdot4(r-s);\\\frac{(40\cdot3^t+4t+2)p^{\prime}-21\cdot3^t-(42\cdot3^t+1)r+p}{8},\frac{(42\cdot3^t+1)r+21\cdot3^t}{8}\Bigg)\ge 0.
\end{multlined}
\end{equation}
Now, iterating the same process $(n-2)$ times for all $n\ge 2$, we have
\begin{equation}\label{eq422}
\begin{multlined}
D_{3^t\cdot2^np^{\prime},3^t\cdot2^ns}\Bigg(L+3^t\cdot2^{n-2}+3^t\cdot2^{n-1}(r-s),L-3^t\cdot2^{n-2}-3^t\cdot2^{n-1}(r-s);\\\frac{\left(\frac{2^{2n+1}-8}{3}\cdot3^t+4t+2\right)p^{\prime}-\frac{2^{2n}-1}{3}\cdot3^t-\left(\frac{2^{2n+1}-2}{3}\cdot3^t+1\right)r+p}{2^n},\\\frac{\left(\frac{2^{2n+1}-2}{3}\cdot3^t+1\right)r+\frac{2^{2n}-1}{3}\cdot3^t}{2^n}\Bigg)\ge 0,    
\end{multlined}    
\end{equation}
which completes the proof of Theorem \ref{thm27}.\qed\\

\section{Acknowledgments}
The authors would like to thank the organizers of the Legacy of Ramanujan $2024$ conference held at Penn State University from June $6$-$9$, $2024$ for their kind invitation where these results were discussed. The authors would also like to thank George E. Andrews and Ali K. Uncu for their kind interest.\\

\section{Declaration}
\subsection*{Ethical Approval}
Not applicable.\\
\subsection*{Funding}
Not applicable.\\

\bibliographystyle{amsplain}


\end{document}